\documentclass[12pt]{amsart}
\usepackage[paperwidth = 8.5in, paperheight = 11in, inner = 1in, outer = 1in, top = 1in, bottom = 1in]{geometry}
\usepackage[pdfborder={0 0 0}]{hyperref}
\hypersetup{
    colorlinks=true,
    linkcolor=blue,
}
\usepackage{graphicx}
\usepackage{amsmath, amsthm, amssymb,amsfonts, mathrsfs}
\usepackage[english]{babel}
\usepackage[utf8]{inputenc}
\usepackage{mathtools}
\usepackage{color}
\usepackage{url}
\newcommand{\FF}{{\mathbb{F}}}

\newcommand{\PP}{{\mathbb{P}}}

  \newcommand{\E}{{\mathcal{E}}}
  \newcommand{\F}{{\mathcal{F}}}
  \newcommand{\G}{{\mathcal{G}}}
\renewcommand{\H}{{\mathcal{H}}}

\renewcommand{\L}{{\mathcal{L}}}
  
\renewcommand{\O}{{\mathcal{O}}}

\title{Computing higher direct images in Macaulay2}
\author{Sasha Zotine}

\date{\today}

\newcommand{\Addresses}{{
  \bigskip
  \footnotesize

  \noindent Sasha Zotine, \textsc{Department of Mathematics \& Statistics,
McMaster University, Hamilton, Ontario, Canada, L8S 4K1}\par\nopagebreak
  \noindent \texttt{E-mail:} \url{zotinea@mcmaster.ca}
}}

\begin{document}

\begin{abstract}
This article highlights the \texttt{ToricHigherDirectImages} package in \textit{Macaulay2}. The central feature is a method for computing (higher) direct images of line bundles under surjective toric morphisms.
\end{abstract}

\maketitle

The package \texttt{ToricHigherDirectImages}, included in version 1.25.05 of \textit{Macaulay2}, is meant to implement the algorithms described in \cite{RothZotineHDI}, wherein Roth and the author provide a method for computing higher direct images of line bundles under toric fibrations. Previously the state-of-the-art for computing higher direct images is in \cite[Corollary~1.3]{EisenbudSchreyer}, which was implemented via the method \texttt{directImageComplex} in the \texttt{BGG} package. However, this method only works for maps between a product of projective spaces. This new package instead allows for maps between any complete toric varieties satisfying a mild assumption. In addition, the package provides an independent algorithm for computing pushforwards of any sheaf under toric Frobenius maps based on \cite[Theorem~2]{AchingerFormula}.

\section*{Higher Direct Images of Toric Maps}

Given a morphism $f \colon X \rightarrow Y$ between varieties the pushforward functor $f_*$ is covariant and left-exact, and hence admits a right derived functor $R^i f_*$. Given a quasi-coherent sheaf $\F$, we may locally describe $R^i f_* \F$ as being given by $R^i f_* \F(U) = H^i(f^{-1}(U), \F)$ for affine open $U \subset Y$, realizing the higher direct images as a generalization of sheaf cohomology to the relative setting; see \cite[Proposition~III.8]{Hartshorne1977}.

Our setting is when $X$ and $Y$ are both smooth complete toric varieties and the map $f$ is a \textit{toric fibration}, i.e. restricts to a group homomorphism on the tori, is equivariant with respect to this action, and satisfies $f_* \O_X = \O_Y$. In this case, each of the varieties, the map, and sheaves on the varieties can easily be encoded in \textit{Macaulay2}, which is done via the \texttt{NormalToricVarieties} package. In particular, the toric varieties $X$ and $Y$ are described using their fans, the map $f$ corresponds to a linear map between the fans which preserves cone inclusions (\cite[Theorem~3.3.4]{CLS2011}), and coherent sheaves on $X$ correspond to graded modules over the Cox ring of $X$.

Here is an example of the central method. Set $Y = \PP^1$ and let $X = \FF_1 = \PP \bigl( \O_Y \oplus \O_Y(1) \bigr)$ be the first Hirzebruch surface. Being a projective bundle, there is a natural projection map $\pi \colon X \rightarrow Y$. The (higher) direct images of a line bundle on a projective bundle are classically known to be vector bundles themselves, which then means they split over the projective line; see \cite[Exercise III.8.4]{Hartshorne1977}. The package recovers this behavior.

\begin{footnotesize}
\begin{verbatim}
Macaulay2, version 1.24.05
with packages: ConwayPolynomials, Elimination, IntegralClosure, InverseSystems, 
Isomorphism, LLLBases, MinimalPrimes, OnlineLookup, PrimaryDecomposition, ReesAlgebra,
Saturation, TangentCone, Truncations,Varieties
i1 : needsPackage "ToricHigherDirectImages"
o1 = ToricHigherDirectImages
o1 : Package
i2 : FF1 = hirzebruchSurface 1; PP1 = toricProjectiveSpace 1; 
     phi = map(PP1,FF1,matrix{{1,0}});
o4 : ToricMap PP1 <--- FF1
i5 : L = OO toricDivisor({0,-1,-2,-3}, FF1)
          1
o5 = OO     (-1, -4)
       FF1
o5 : coherent sheaf on FF1, free of rank 1
i6 : prune phi_*^0 L
o6 = 0
o6 : coherent sheaf on PP1, free of rank 0
i7 : prune phi_*^1 L
          1             1             1
o7 = OO    (-2) ++ OO    (-3) ++ OO    (-4)
       PP1           PP1           PP1
o7 : coherent sheaf on PP1, free of rank 3
\end{verbatim}
\end{footnotesize}

Higher direct images need not always split as a direct sum of line bundles. For instance, the first Hirzebruch surface $\FF_1$ may also be interpreted as the blowup of a torus fixed point on $\PP^2$, and the pushforward of $\O_{\FF_1}(-2E)$, where $E$ is the exceptional divisor, gives the ideal sheaf of the blown up point.

\begin{footnotesize}
\begin{verbatim}
i8 : PP2 = toricProjectiveSpace 2; psi = map(PP2,FF1,matrix{{0,-1},{1,0}});
o9 : ToricMap PP2 <--- FF1
i10 : E = {0,1,0,0};
i11 : prune HDI(psi,0,-2*E)
o11 = cokernel {2} | x_0  0    |
               {2} | -x_2 x_0  |
               {2} | 0    -x_2 |
                                                 3
o11 : QQ[x ..x ]-module, quotient of (QQ[x ..x ])
          0   2                           0   2
\end{verbatim}
\end{footnotesize}

There is also behavior involving both torsion and free summands; here is \cite[Example~1.8]{RothZotineHDI}.

\begin{footnotesize}
\begin{verbatim}
i12 : X = toricBlowup({1,5}, FF1 ** PP1);
i13 : theta = map(FF1, X, matrix{{1,0,0},{0,1,0}});
o13 : ToricMap FF1 <--- X
i14 : D = {0,0,0,0,-2,0,-2};
i15 : prune theta_*^1 D
o15 = cokernel | 0   0     |
               | 0   x_1^2 |
               | x_1 0     |
                                                 3
o15 : QQ[x ..x ]-module, quotient of (QQ[x ..x ])
          0   3                           0   3
\end{verbatim}
\end{footnotesize}

Higher direct images enjoy a relative version of Serre duality, sometimes known as Grothendieck-Serre duality or coherent duality. For a proper morphism $f \colon Z \rightarrow Y$ and coherent sheaves $\F$ on $Z$ and $\G$ on $Y$, Grothendieck duality relates the derived Hom functor ($\E xt$) and higher direct images via the equation 
\begin{equation}
\label{eq:duality}
    \E xt_Y(Rf_* \F, \G) = R f_* \E xt_Z(\F, f^! \G),
\end{equation} 
see \cite[Ideal~Theorem~(c), p.~7]{Hartshorne1966residues}. While this statement holds for complexes of sheaves and our package does not yet implement functionality for working with complexes, we may nevertheless witness this duality term-by-term. We consider a proper surjection $f \colon Z \rightarrow \PP^4$ from a smooth Fano sixfold to $\PP^4$. We set $\F$ to be a line bundle whose higher direct images are only supported in degree two.

\begin{footnotesize}
\begin{verbatim}
i16 : Z = smoothFanoToricVariety(6,174); f = (nefRayContractions Z)_2; PP4 = target f;
o17 : ToricMap normalToricVariety ({{-1,0,0,0},{0,-1,0,0},{0,0,-1,1}, {0,0,0,-1},
          {1,1,1,0}},{{0,1,2,3},{0,1,2,4},{0,1,3,4},{0,2,3,4},{1,2,3,4}}) <--- Z
i19 : F = OO_Z(0,-4,1,0);
i20 : assert(f_*^0 F == 0 and f_*^1 F == 0);
i21 : prune f_*^2 F
o21 = cokernel {2} | x_1 0   |
               {2} | 0   x_1 |
               {4} | 0   0   |
                                              2             1
o21 : coherent sheaf on PP4, quotient of OO    (-2) ++ OO    (-4)
                                           PP4           PP4
\end{verbatim}
\end{footnotesize}

Set $\G = K_Y[\dim Y]$ to be the canonical line bundle on $Y$, placed in degree $-\dim Y = -4$ as a complex. We compute the right side of \eqref{eq:duality}. We claim that $f^! \G = K_Z[\dim Z]$. Indeed, for a map of any scheme to a point $g \colon X \rightarrow \{\text{pt}\}$, we have that $g^! \O_{\{\text{pt}\}} = K_X[\dim X]$ using \cite[Definition of $f^\#$]{Hartshorne1966residues} and \cite[Theorem~8.7~(4)]{Hartshorne1966residues}. Since shrieks are functorial, we may take $X = \PP^4$ and get
\begin{equation*}
    f^! \G = f^! (g^! \O_{\{\text{pt}\}}) = (f \circ g)^! \O_{\{\text{pt}\}} = K_Z[\dim Z].
\end{equation*}
Now for any locally free sheaf $\E$ on $Z$, we have that $\E xt_Z^0(\E, K_Z) = \H om(\E, K_Z) = \E^* \otimes K_Z$. Combining these facts together, \eqref{eq:duality} turns into
\begin{equation*}
    \E xt_{\PP^4}(Rf_* \F, K_{\PP^4})[-2] = Rf_*(\F \otimes K_Z),
\end{equation*}
which we may verify computationally.

\begin{footnotesize}
\begin{verbatim}
i22 : KZ = OO toricDivisor Z; KPP4 = OO toricDivisor PP4;
i24 : assert(sheaf_PP4 Ext^0(module prune f_*^2 F, module KPP4) == f_*^0 (KZ ** dual F));
i25 : assert(sheaf_PP4 Ext^1(module prune f_*^2 F, module KPP4) == f_*^1 (KZ ** dual F));
\end{verbatim}
\end{footnotesize}

In the example above, we make use of the method \texttt{nefRayContractions} in the package for easily producing surjective morphisms. As far as the author is aware, this is not the first package to implement such a construction, which comes from a procedure in the toric minimal model program; see \cite[Chapters 14 and 15]{CLS2011} for details.

Another method which is exported is \texttt{computeEigencharacters}, which packages the data output by \cite[Algorithm~A.8]{RothZotineHDI} into a hash table. More precisely, for a toric morphism $f \colon X \rightarrow Y$ and line bundle $\L$ on $X$, the higher direct image $R^i f_* \L$ split into eigensheaves indexed by characters over the kernel torus $T_K = \ker f|_{T_X}$. This method gives a hash table whose keys are these characters and values are the divisors on $X$ and $Y$ used for computing higher direct images. For example, we reproduce the computation in \cite[Example~B.3]{RothZotineHDI}.
\begin{footnotesize}
\begin{verbatim}
i26 : keys computeEigencharacters(theta,1,D)
o26 = {| 0 |, 0, | 0  |}
       | 0 |     | 0  |
       | 1 |     | -1 |
o26 : List
\end{verbatim}
\end{footnotesize}

\section*{Toric Frobenius Pushforwards}

Given any toric variety $X$, the \textit{$p$th toric Frobenius map} $F_p \colon X \rightarrow X$ is a surjective toric morphism induced by the map $T_X \rightarrow T_X$ sending a torus element $(t_1, t_2, \ldots, t_n)$ to $(t_1^p, t_2^p, \ldots, t_n^p)$. This map is defined irrespective of the characteristic of the ground field, though was classically introduced in the context of finite characteristics. Being a finite map, the toric Frobenius maps have no higher direct images. Their direct images, on the other hand, are of some recent interest; see the Bondal-Thomsen collection in \cite[Section~2]{BBBHEFGHH2024King}. Thomsen showed in \cite[Theorem~1.1]{ThomsenFrobenius} that the Frobenius pushforward of a line bundle splits, and the summands which appear can be combinatorially described; see Bondal's article for details \cite{BondalFormula}. The method \texttt{frobeniusDirectImage} computes this splitting using an adapted algorithm from \cite{AchingerFormula}. Moreover, we also introduce functoriality, allowing for the pushforward of maps between line bundles. Consequently, we are able to compute the pushforward of any module by resolving it by line bundles and pushing forward the complex. For instance, on the first Hirzebruch surface, we can compute the Frobenius pushforward of the cotangent sheaf.
\begin{footnotesize}
\begin{verbatim}
i27 : Omega = module cotangentSheaf FF1
o27 = cokernel {2, 0} | x_1x_3 |
               {0, 2} | x_0    |
               {0, 2} | -x_2   |
                                                 3
o27 : QQ[x ..x ]-module, quotient of (QQ[x ..x ])
          0   3                           0   3
i28 : prune frobeniusDirectImage(2, Omega)
o28 = cokernel {1, 0} | 0      |
               {2, 0} | x_1x_3 |
               {0, 1} | 0      |
               {0, 2} | x_0    |
               {1, 1} | 0      |
               {0, 1} | 0      |
               {0, 2} | 0      |
               {1, 1} | 0      |
               {0, 2} | -x_2   |
                                                 9
o28 : QQ[x ..x ]-module, quotient of (QQ[x ..x ])
          0   3                           0   3
\end{verbatim}
\end{footnotesize}
As mentioned, this presentation is obtained from the free resolution of the cotangent bundle.
\begin{footnotesize}
\begin{verbatim}
i29 : (frobeniusDirectImage(2, freeResolution Omega)).dd
                      12                                                    4
o29 = 0 : (QQ[x ..x ])   <------------------------------------- (QQ[x ..x ])  : 1
               0   3        {1, 0} | 0   0      0    x_1x_3 |        0   3
                            {1, 1} | 0   0      1    0      |
                            {2, 0} | 0   x_1x_3 0    0      |
                            {1, 1} | 1   0      0    0      |
                            {0, 1} | x_0 0      0    0      |
                            {0, 2} | 0   x_0    0    0      |
                            {1, 1} | 0   0      1    0      |
                            {0, 2} | 0   0      0    1      |
                            {0, 1} | 0   0      -x_2 0      |
                            {0, 2} | 0   0      0    -1     |
                            {1, 1} | -1  0      0    0      |
                            {0, 2} | 0   -x_2   0    0      |
o29 : ComplexMap
\end{verbatim}
\end{footnotesize}

\section*{Acknowledgements}

I am very grateful to the Fields Institute for hosting me during the writing of this article. My thanks to Mahrud Sayrafi and Gregory G.~Smith for helpful conversations on writing the package. I am also very thankful to Mike Roth, who provided me with the nice example of Grothendieck duality. This research was partially funded by the Natural Sciences and Engineering Research Council (NSERC).

\bibliographystyle{plain}
\bibliography{refs}

\begin{thebibliography}{1}

\bibitem{AchingerFormula}
Piotr Achinger.
\newblock A characterization of toric varieties in characteristic {$p$}.
\newblock {\em Int. Math. Res. Not. IMRN}, 16:6879--6892, 2015.

\bibitem{BBBHEFGHH2024King}
Matthew~R. Ballard, Christine Berkesch, Michael~K. Brown, Lauren~Cranton Heller, Daniel Erman, David Favero, Sheel Ganatra, Andrew Hanlon, and Jesse Huang.
\newblock King's conjecture and birational geometry.
\newblock {\em arXiv}, 2024.
\newblock Available at \href{https://arxiv.org/abs/2501.00130}{\texttt{arXiv:2501.00130}}.

\bibitem{BondalFormula}
Alexey Bondal.
\newblock Derived categories of toric varieties.
\newblock Oberwolfach {Report} 5/2006, {No}. 1, 284-286, 2006.

\bibitem{CLS2011}
David~A. Cox, John~B. Little, and Henry~K. Schenck.
\newblock {\em Toric Varieties}.
\newblock Graduate studies in mathematics. American Mathematical Society, 2011.

\bibitem{EisenbudSchreyer}
David Eisenbud and Frank-Olaf Schreyer.
\newblock Relative {B}eilinson monad and direct image for families of coherent sheaves.
\newblock {\em Trans. Amer. Math. Soc.}, 360(10):5367--5396, 2008.

\bibitem{Hartshorne1966residues}
Robin Hartshorne.
\newblock {\em Residues and duality}, volume No. 20 of {\em Lecture Notes in Mathematics}.
\newblock Springer-Verlag, Berlin-New York, 1966.
\newblock Lecture notes of a seminar on the work of A. Grothendieck, given at Harvard 1963/64, With an appendix by P. Deligne.

\bibitem{Hartshorne1977}
Robin Hartshorne.
\newblock {\em Algebraic Geometry}.
\newblock Graduate Texts in Mathematics. Springer, 1977.

\bibitem{RothZotineHDI}
Mike Roth and Sasha Zotine.
\newblock Reduced \v{C}ech complexes and computing higher direct images under toric maps.
\newblock {\em arXiv}, 2025.
\newblock Available at \href{https://arxiv.org/abs/2504.12903}{\texttt{arXiv:2504.12903}}.

\bibitem{ThomsenFrobenius}
Jesper~Funch Thomsen.
\newblock Frobenius direct images of line bundles on toric varieties.
\newblock {\em J. Algebra}, 226(2):865--874, 2000.

\end{thebibliography}

\Addresses

\end{document}